\documentclass{amsart}
\usepackage{amsfonts,amssymb,amsmath,amsthm}
\usepackage{url}
\usepackage{enumerate}
\newtheorem{thm}{Theorem}[section]
\newtheorem{cor}[thm]{Corollary}
\newtheorem{lem}[thm]{Lemma}
\newtheorem{prop}[thm]{Proposition}
\theoremstyle{definition}

\theoremstyle{remark}

\newtheorem{exe}[thm]{\bf Example}
\numberwithin{equation}{section}
\begin{document}
\title[Running Head]{Hypercyclic abelian affine groups }

\author{Adlene Ayadi}

\address{Adlene Ayadi,  Department of Mathematics,
Faculty of Sciences of Gafsa,University of Gafsa, Gafsa, Tunisia}

\email{adleneso@yahoo.fr}

\thanks{This work is supported by the research unit: syst\`emes dynamiques et combinatoire: 99UR15-15}



\begin{abstract}

In this paper, we give a characterization of hypercyclic abelian affine group $\mathcal{G}$. If $\mathcal{G}$
is finitely generated, this characterization is explicit. We prove in particular
 that no abelian group generated by $n$ affine maps on $\mathbb{C}^{n}$ has a dense orbit.
\end{abstract}

\maketitle

\section{ Introduction}

Let $M_{n}(\mathbb{C})$ be the set of all square matrices of order
$n\geq 1$  with entries in  $\mathbb{C}$ and $GL(n, \ \mathbb{C})$
be the group of all invertible matrices of $M_{n}(\mathbb{C})$. A
map $f: \ \mathbb{C}^{n}\longrightarrow \mathbb{C}^{n}$  is called
an affine map if there exist $A\in M_{n}(\mathbb{C})$  and  $a\in
\mathbb{C}^{n}$  such that  $f(x)= Ax+a$,  $x\in \mathbb{C}^{n}$.
We denote  $f= (A,a)$, we call  $A$ the \textit{linear part} of
$f$. Denote by  $MA(n, \ \mathbb{C})$  the set of all affine maps
and $GA(n, \ \mathbb{C})$ the set of all invertible affine maps of
$MA(n, \mathbb{C})$. $MA(n, \mathbb{C})$ is a vector space and for
composition of maps,  $GA(n, \ \mathbb{C})$ is a group.

Let  $\mathcal{G}$  be an abelian affine subgroup of $GA(n, \
\mathbb{C})$. For a vector  $v\in \mathbb{C}^{n}$,  we consider
the orbit of  $\mathcal{G}$  through  $v$: $ \mathcal{G}(v) =
\{f(v): \ f\in \mathcal{G}\} \subset \mathbb{\mathbb{C}}^{n}$. A
subset $E\subset \mathbb{C}^{n}$  is called \emph{$
\mathcal{G}$-invariant} if  $f(E)\subset E$  for any  $f\in
\mathcal{G}$;  that is  $E$  is a union of orbits. Before stating
our main results, we introduce the following notions:

A subset  $\mathcal{H}$  of  $\mathbb{C}^{n}$  is called {\it an
affine subspace} of  $\mathbb{C}^{n}$  if there exist a vector
subspace  $H$  of  $\mathbb{C}^{n}$  and  $a\in \mathbb{C}^{n}$
 such that  $\mathcal{H} = H + a$. For  $a\in \mathbb{C}^{n}$,
denote by  $T_{a}: \ \mathbb{C}^{n}\longrightarrow
\mathbb{C}^{n};$
 $x\longmapsto x+a$  the translation map by vector  $a$, so
$\mathcal{H} = T_{a}(H)$.  We say that  $\mathcal{H}$  has
dimension $p$  ($0\leq p\leq n$), denoted dim($\mathcal{H})=p$, if
$H$  has dimension  $p$.
\medskip

Denote by $\overline{A}$ the closure of a  subset  $A \subset
\mathbb{C}^{n}$. A subset  $E$  of  $\mathbb{C}^{n}$   is called
\emph{a minimal set} of  $\mathcal{G}$  if  $E$  is closed in  $
\mathbb{C}^{n}$,  non empty,  $ \mathcal{G}$-invariant and has no
proper subset with these properties. It is equivalent to say that
$E$  is a $\mathcal{G}$-invariant set such that every orbit
contained in  $E$ is dense in it.
The group $\mathcal{G}$ is
called \textit{hypercyclic}
if there exists a vector $v\in {\mathbb{C}}^{n}$ such that $\mathcal{G}(v)$
is dense  in ${\mathbb{C}}^{n}$. For an account of results and bibliography on hypercyclicity, we refer to
the book \cite{bm} by Bayart and Matheron.
\\
\\
Define the map  $$\Phi\ :\ GA(n,\ \mathbb{C})\ \longrightarrow\
\Phi(GA(n,\mathbb{C}))\subset GL(n+1,\ \mathbb{C})$$ \ \ \ \ \ \
\ \ \ \ \ \ \  \ \ \ \ \ \ \ \ \ \ \ \ \ \  \ \ \ \ \ \ \  \ \ \ \
\ \ \  \ $f =(A,a) \ \longmapsto\ \begin{bmatrix}
                     1  & 0 \\
                     a & A
                 \end{bmatrix}$
\\
We have the following composition formula
 $$\begin{bmatrix}
                     1  & 0 \\
                     a & A
                 \end{bmatrix}\begin{bmatrix}
                     1  & 0 \\
                     b & B
                 \end{bmatrix} = \begin{bmatrix}
                     1  & 0 \\
                     Ab+a & AB
                 \end{bmatrix}.$$
 Then  $\Phi$  is a homomorphism of groups.
                 \
                 \\
                 \\
 Let $\mathcal{G}$  be an abelian affine subgroup of
$GA(n, \ \mathbb{C})$. Then  $\Phi(\mathcal{G})$  is an abelian
subgroup of  $GL(n+1, \mathbb{C})$.
\\
\\
Denote by:
\\
\textbullet \ $\mathbb{C}^{*}= \mathbb{C}\backslash\{0\}$ and
$\mathbb{N}_{0}= \mathbb{N}\backslash\{0\}$.
\
\\
Let $n\in\mathbb{N}_{0}$  be fixed. For each $m=1,2,\dots, n+1$,
denote by:\
\\
\textbullet \;  $\mathcal{B}_{0} = (e_{1},\dots,e_{n+1})$ the
canonical basis of $\mathbb{C}^{n+1}$  and   $I_{n+1}$  the
identity matrix of $GL(n+1,\mathbb{C})$.
\
\\
\\
\textbullet \; $\mathbb{T}_{m}(\mathbb{C})$ the set of matrices
over $\mathbb{C}$ of the form $$\begin{bmatrix}
  \mu &  &   & 0 \\
  a_{2,1} & \mu &   &  \\
  \vdots &  \ddots & \ddots &  \\
  a_{m,1} & \dots & a_{m,m-1} & \mu
\end{bmatrix} \qquad (1)$$
\\
\\
\textbullet \; $\mathbb{T}_{m}^{\ast}(\mathbb{C})$  the group of
matrices of the form $(1)$ with $\mu\neq 0$.
\
\\
\\
Let  $r\in \mathbb{N}$ and $\eta =(n_{1},\dots,n_{r})$
 be a sequence of positive integers such that $n_{1}+\dots + n_{r}=n+1. $
In particular, $r\leq n+1$.
\
\\
Write
\
\\
\textbullet \; \ $\mathcal{K}_{\eta,r}(\mathbb{C}): =
\mathbb{T}_{n_{1}}(\mathbb{C})\oplus\dots \oplus
\mathbb{T}_{n_{r}}(\mathbb{C}).$ In particular if  $r=1$, then
$\mathcal{K}_{\eta,1}(\mathbb{C}) = \mathbb{T}_{n+1}(\mathbb{C})$
and  $\eta=(n+1)$.
\
\\
\textbullet \; $\mathcal{K}^{*}_{\eta,r}(\mathbb{C}): =
\mathcal{K}_{\eta,r}(\mathbb{C})\cap \textrm{GL}(n+1, \
\mathbb{C})$.
\
\\
\\
Define the map  $$\Psi\ :\ MA(n,\ \mathbb{C})\ \longrightarrow\
\Psi(MA(n,\mathbb{C}))\subset M_{n+1}(\mathbb{C})$$ \ \ \ \ \ \  \
\ \ \ \ \ \  \ \ \ \ \ \ \ \ \ \ \ \ \ \  \ \ \ \ \ \ \  \ \ \ \ \
\ \  \ $f =(A,a) \ \longmapsto\ \begin{bmatrix}
                     0  & 0 \\
                     a & A
                 \end{bmatrix}$
\\
We have   $\Psi$  is an isomorphism.
\
\\
\textbullet \; $\mathcal{F}_{n+1}=\Psi(MA(n,\mathbb{C}))$.
\
\\
\textbullet \; $\textrm{exp} :\ \mathbb{M}_{n+1}(\mathbb{C})
\longrightarrow\textrm{GL}(n+1, \mathbb{C})$  is the matrix
exponential map; set $\textrm{exp}(M) = e^{M}$.\\

There always exists a $P\in \Phi(\textrm{GA}(n, \mathbb{C}))$ and
a partition $\eta$ of $n+1$ such that $G'=P^{-1}GP\subset
\mathcal{K}^{*}_{\eta,r}(\mathbb{C})\cap\Phi(GA(n,\mathbb{C}))$
(see Proposition ~\ref{p:2}). For such a choice of matrix $P$,  we
let
\\
\\
\textbullet\  $\mathrm{g} = \textrm{exp}^{-1}(G)\cap \left(
P(\mathcal{K}_{\eta,r}(\mathbb{C}))P^{-1}\right)\cap
\mathcal{F}_{n+1}$. If $G\subset
\mathcal{K}^{*}_{\eta,r}(\mathbb{C})$,  we have $\mathrm{g} =
\textrm{exp}^{-1}(G)\cap \mathcal{K}_{\eta,r}(\mathbb{C})\cap
\mathcal{F}_{n+1}$.
\
\\
\textbullet\  $\mathrm{g}_{u} = \{Bu: \ B\in \mathrm{g}\}, \
u\in\mathbb{C}^{n+1}.$
\
\\
\textbullet\  $\mathfrak{g} = \Psi^{-1}(\mathrm{g})$.
\\
\textbullet\  $\mathfrak{g}_{v} =\{f(v),\ \ f\in \mathfrak{g}\}, \
v\in\mathbb{C}^{n}.$
\
\\
\textbullet \; $u_{0} = [e_{1,1},\dots,e_{r,1}]^{T}\in
\mathbb{C}^{n+1}$ where
 $e_{k,1} = [1,0,\dots,
0]^{T}\in \mathbb{C}^{n_{k}}$, for $k=1,\dots, r$. One has
$u_{0}\in\{1\}\times\mathbb{C}^{n}$.
\\
\textbullet \;
$p_{2}:\mathbb{C}^{n+1}\longrightarrow\mathbb{C}^{n}$ the
projection defined by
$p_{2}(x_{1},\dots,x_{n+1})=(x_{2},\dots,x_{n+1})$.
 \\
\textbullet \; $v_{0} = Pu_{0}$. As
$P\in\Phi(GA(n, \mathbb{C}))$,
$v_{0}\in \{1\}\times\mathbb{C}^{n}$.
 \\
 \textbullet \; $w_{0} = p_{2}v_{0}\in\mathbb{C}^{n}$.
 \\
 \textbullet \; $e^{(k)} = [e^{(k)}_{1},\dots,  e^{(k)}_{r}]^{T}\in \mathbb{C}^{n+1}$ where $$e^{(k)}_{j} = \left\{\begin{array}{c}
  0\in \mathbb{C}^{n_{j}}\ \ \mathrm{if}\ \ j\neq k \\
  e_{k,1}\ \ \ \ \ \ \ \ \mathrm{if}\ \ j = k  \\
\end{array}\right. \ \ \ \ \ \  for \ \ every\ \ 1\leq j,k\leq r.$$

\bigskip

For groups of affine maps on  $\mathbb{K}^{n}$
($\mathbb{K}=\mathbb{R}$ or $\mathbb{C}$), their dynamics were
recently initiated for some classes in different point of view,
(see for instance, \cite{Ja}, \cite{Ku}, \cite{be},\cite{Bu}). The
purpose here is to give analogous results of that
 theorem for linear abelian subgroup of $GL(n, \mathbb{C})$ proved in \cite{aAh-M05} (see Proposition ~\ref{C0C:1}).
  Our main results are the following:
\medskip

\begin{thm}\label{T:1} Let  $\mathcal{G}$ be an abelian subgroup of $GA(n,
\mathbb{C})$. The following are equivalent:
\
\\
(i) $\mathcal{G}$ is hypercyclic.
\\
(ii) The orbit $\mathcal{G}(w_{0})$ is dense in $\mathbb{C}^{n}$
\\
(iii) $\mathfrak{g}_{w_{0}}$ is an additive subgroup dense in
$\mathbb{C}^{n}$
\end{thm}
\bigskip

For a \textit{finitely generated} abelian subgroup $\mathcal{G}\subset
\textrm{GA}(n, \mathbb{R})$, let introduce the following property.
Consider the following rank condition on a collection of affine
maps $f_{1},\dots,f_{p}\in
GA(n, \mathbb{C})$, where  $f'_{1},\dots,f'_{p}\in
\mathfrak{g}$ such that $e^{\Psi(f'_{k})} = \Phi(f_{k})$,
$k=1,\dots,p$.
\\
\\
We say that $f_{1},\dots,f_{p}$ satisfy \emph{property
$\mathcal{D}$} if  for every $(s_{1},\dots,s_{p};\
t_{1},\dots,t_{r})\in \mathbb{Z}^{p+r}\backslash\{0\}$: \
 $$rank
\left[\begin{array}{cccccc}
   Re(f'_{1}(w_{0})) & \dots& Re(f'_{p}(w_{0}))& 0 &\dots &0\\
      Im(f'_{1}(w_{0})) & \dots& Im(f'_{p}(w_{0}))& 2\pi p_{2}\circ Pe^{(1)}&\dots &2\pi p_{2}\circ Pe^{(r)} \\
      s_{1} & \dots& s_{p}& t_{1}&\dots& t_{r} \\
 \end{array}\right]=2n+1.$$\

For a vector $v\in\mathbb{C}^{n}$, we write $v = Re(v)+ iIm(v)$
where $Re(v)$ and $Im(v)\in \mathbb{R}^{n}$. In this case, the Theorem can
be stated as follows:
\bigskip

\begin{thm}\label{T:2} Let  $\mathcal{G}$  be an abelian subgroup of  $GA(n, \mathbb{C})$  generated by
 $f_{1}, \dots,f_{p}$
and let $f'_{1},\dots,f'_{p}\in\mathfrak{g}$ such that $e^{\Psi(f'_{1})} =
\Phi(f_{1}),\dots, e^{\Psi(f'_{p})} = \Phi(f_{p})$. Then the
following are equivalent:
\
\\
(i)  $\mathcal{G}$ is hypercyclic.
\\
(ii)  the maps $f_{1},\dots,f_{p}$ satisfy property $\mathcal{D}$
\\
(iii)
$\mathfrak{g}_{w_{0}}=  \underset{k=1}{\overset{p}{\sum}}\mathbb{Z}f'_{k}(w_{0})+2i\pi\underset{k=1}{\overset{r}{\sum}}\mathbb{Z}(p_{2}\circ
Pe^{(k)})$ is an additive group dense in $\mathbb{C}^{n}$.
\end{thm}
\medskip

\begin{cor}\label{C:9} If $\mathcal{G}$ is of finite type $p$ with $p \leq 2n - r+1$, then it has no
dense orbit.
\end{cor}
\medskip

\begin{cor}\label{C:10} If $\mathcal{G}$ is of finite type $p$ with $p\leq n$, then it has no dense
orbit.
\end{cor}
\medskip

\section{Notations and Lemmas}

Denote by  $\mathcal{L}_{\mathcal{G}}$
 the set of the linear parts of all elements of  $\mathcal{G}$ and $vect(F)$ is the vector space generated by a subset $F\subset \mathbb{C}^{n+1}$.
In the following, denote by $I_{m}$ the identity matrix of $GL(m, \mathbb{C})$, for any $m\in \mathbb{N}_{0}$.

\begin{prop}\label{p:1}$($\cite{aAh-M05}, Proposition 2.3$)$ Let  $L$ be an abelian subgroup of  $GL(n,\mathbb{C})$. Then
there exists $P\in GL(n,\mathbb{C})$ such that\ $P^{-1}LP$ is a
subgroup of  $\mathcal{K}_{\eta',r'}(\mathbb{C})$, for some $r'\leq n$ and $\eta'=(n'_{1},\dots,n'_{r'})\in\mathbb{N}_{0}^{r'}$.
\end{prop}

\medskip

\begin{prop}\label{p:2} Let  $\mathcal{G}$ be an abelian subgroup of
$GA(n,\mathbb{C})$ and  $G=\Phi(\mathcal{G})$. Then
there exists $P\in \Phi(GA(n,\mathbb{C}))$ such that $P^{-1}GP$ is a subgroup of $\mathcal{K}^{*}_{\eta,r}(\mathbb{C})\cap\Phi(GA(n, \mathbb{C}))$,
 for some $r\leq n+1$ and $\eta=(n_{1},\dots,n_{r})\in\mathbb{N}_{0}^{r}$. In particular, $Pu_{0}\in \{1\}\times \mathbb{C}^{n}$.
\end{prop}
\medskip

\begin{proof} We have $\mathcal{L}_{\mathcal{G}}$ is an abelian subgroup of $GL(n, \mathbb{C})$. By Proposition ~\ref{p:1}, there exists
$Q\in GL(n,\mathbb{C})$ such that\ $Q^{-1}\mathcal{L}_{G}Q$ is a
subgroup of  $\mathcal{K}^{*}_{\eta',r'}(\mathbb{C})$ for some $r'\leq n$ and $\eta'=(n'_{1},\dots,n'_{r'})\in\mathbb{N}_{0}^{r'}$ such that $n'_{1}+\dots,n'_{r'}=n$.
For every $A\in\mathcal{L}_{\mathcal{G}}$,  $Q^{-1}AQ=\mathrm{diag}(A_{1},\dots,A_{r'})$ with $A_{k}\in \mathbb{T}^{*}_{n'_{k}}$ and $\mu_{A_{k}}$ is the only eigenvalue
of $A_{k}$, $k=1,\dots,r'$. Let $J=\{k\in\{1,\dots,r'\},\ \ \mu_{A_{k}}=1, \ \forall\ A\in \mathcal{L}_{\mathcal{G}}\}$. There are two cases:
\\
\\
- {\it Case1:} Suppose that $J=\{k_{1},\dots,k_{s}\}$ for some $s\leq r'$.  We can take $J=\{1,\dots,s\}$, otherwise, we replace $P_{1}$ by $RP_{1}$ for some permutation
 matrix $R$  of $GL(n, \mathbb{C})$. Let  $P_{1}=\mathrm{diag}(1,Q)$, so $P_{1}\in\Phi(GA(n,\mathbb{C}))$ and $$P_{1}^{-1}\Phi(f)P_{1}=\left[\begin{array}{cc}
                                                        1 & 0 \\
                                                        Q^{-1}a & Q^{-1}AQ
                                                      \end{array}
 \right].$$
 Write  $E=\mathrm{vect}(\mathcal{C}_{1},\dots,\mathcal{C}_{s})$ and $H=\mathrm{vect}(\mathcal{C}_{s+1},\dots,\mathcal{C}_{r'})$, so
 $E$ and $H$ are $G$-invariant vector spaces. Moreover, for every $f=(A,a)\in \mathcal{G}$ the restriction $A_{/E}$
  has $1$ as only eigenvalue and we have
$$P_{1}^{-1}\Phi(f)P_{1}=\left[\begin{array}{ccc}
                                                        1&0&0\\
                                                        a_{1}& A_{1} & 0 \\
                                                       a_{2} &0& A_{2}
                                                      \end{array}
 \right], \ \ \ (1)$$ where $A_{1}=A_{/E}\in\mathbb{T}^{*}_{p}(\mathbb{C})$, $A_{2}=A_{/H}\in\mathbb{T}^{*}_{n-p}(\mathbb{C})$, $a_{1}\in \mathbb{C}^{p}$, $a_{2}\in\mathbb{C}^{n-p}$,
 $p=n'_{1}+\dots+n'_{s}$. On the other hand, there exists $f_{0}=(B,b)\in \mathcal{G}$ such that $B_{2}=B_{/H}$ has no eigenvalue equal to $1$,
  so $B_{2}-I_{n-p}$ is invertible. As in (1), write $$P_{1}^{-1}\Phi(f_{0})P_{1}=\left[\begin{array}{ccc}
                                                        1&0&0\\
                                                        b_{1}& B_{1} & 0 \\
                                                       b_{2} &0& B_{2}
                                                      \end{array}
 \right].$$ Let $P_{2}=\left[\begin{array}{ccc}
                                                       1 & 0 & 0 \\
                                                       0 & I_{p} & 0 \\
                                                       b_{2} & 0 & B_{2}-I_{n-p}
                                                     \end{array}
  \right]$ and $P=P_{1}P_{2}^{-1}$, then $P=\left[\begin{array}{cc}
                                                    1 & 0 \\
                                                    d & P_{0}
                                                  \end{array}
  \right]\in\Phi(GA(n, \mathbb{C}))$, where $P_{0}=Q.Q_{1}^{-1}$, $Q_{1}=\left[\begin{array}{cc}
                                                    I_{p} & 0 \\
                                                    0 & B_{2}-I_{n-p}
                                                  \end{array}
  \right]$ and $d=-P_{0}[0,b_{2}]^{T}$.
\
\\
  For every $f=(A,a)\in \mathcal{G}$ we have by $(1)$, \begin{align*}
  P^{-1}\Phi(f)P & =P_{2}P_{1}^{-1}\Phi(f)P_{1}P_{2}^{-1}\\
  \ & =P_{2}\left[\begin{array}{ccc}
                                                        1&0&0\\
                                                        a_{1}& A_{1} & 0 \\
                                                       a_{2} &0& A_{2}
                                                      \end{array}
 \right]P_{2}^{-1}\\
 \ & =\left[\begin{array}{ccc}
 1&0&0\\
 a_{1}& A_{1} & 0 \\
 -(A_{2}-I_{n-p})b_{2}+(B_{2}-I_{n-p})a_{2} &0& A_{2}
  \end{array}
 \right]\ \ \ (2)
\end{align*}

Since $G$ is abelian, so by the equality $P_{1}^{-1}\Phi(f)\Phi(f_{0})P_{1}=P_{1}^{-1}\Phi(f_{0})\Phi(f)P_{1}$, we find
  \ $-(A_{2}-I_{n-p})b_{2}+(B_{2}-I_{n-p})a_{2}=0$.

\
\\
  It follows by (2), that
  $P^{-1}GP$ is a subgroup of $\mathcal{K}^{*}_{\eta,r}(\mathbb{C})\cap \Phi(GA(n, \mathbb{C}))$, where $r=r'-s+1$ and
  $\eta=(p+1,n'_{s+1},\dots,n'_{r'})$.\
  \\
  \\
  - {\it Case2:} Suppose that $J=\emptyset$, and denote by  $Fix(G)=\{x\in\mathbb{C}^{n+1}: Bx=x, \forall B\in G\}$, then $Fix(G)=\mathbb{C}v$
  for some $v=(1, v_{1})$, $v_{1}\in\mathbb{C}^{n}$. For every $f=(A,a)\in\mathcal{G}$, one has $\Phi(f)(1,v_{1})=(1, f(v_{1}))=(1,v_{1})$ so
  $f(v_{1})=Av_{1}+a=v_{1}$.    Let
   $P=\left[\begin{array}{cc}
              1 & 0 \\
              v_{1} & P_{1}
            \end{array}
   \right]\in\Phi(GA(n,\mathbb{C}))$, then  \begin{align*}
   P^{-1}\Phi(f)P & =\left[\begin{array}{cc}
   1 & 0 \\
   -P_{1}^{-1}v_{1}& P^{-1}_{1}
   \end{array}
 \right] \left[\begin{array}{cc}
   1 & 0 \\
   a & A
   \end{array}
 \right] \left[\begin{array}{cc}
   1 & 0 \\
   v_{1}& P_{1}
   \end{array}
 \right]    \\
 \ &\  =  \left[\begin{array}{cc}
 1 & 0 \\
 P_{1}^{-1}(Av_{1}+a-v_{1}) & P_{1}^{-1}AP_{1}
 \end{array}
 \right] \\
 \ &\ =\left[\begin{array}{cc}
 1 & 0 \\
 0 & P_{1}^{-1}AP_{1}
 \end{array}
 \right].
 \end{align*} It follows that
  $P^{-1}GP$ is a subgroup of $\mathcal{K}^{*}_{\eta,r}(\mathbb{C})\cap\Phi(GA(n, \mathbb{C}))$, where $r=r'+1$ and $\eta=(1,n'_{1},\dots,n'_{r'})$.\
\\
Since $u_{0}\in\{1\}\times\mathbb{C}^{n}$ and $P\in\Phi(GA(n, \mathbb{C}))$, so $Pu_{0}\in\{1\}\times\mathbb{C}^{n}$.
\end{proof}
\
\
\\
\\
Denote by:
\
\\
$\bullet$ $G'=P^{-1}GP$.\
\\
$\bullet$ $\mathrm{g}'=exp^{-1}(G')\cap\mathcal{K}_{\eta,r}(\mathbb{C})\cap\mathcal{F}_{n+1}$.\
\\
$\bullet$ $\mathrm{g}_{1}=exp^{-1}(G)\cap\left(P\mathcal{K}_{\eta,r}(\mathbb{C})P^{-1}\right)$
\\
\begin{lem}\label{L:4-}$($\cite{aAh-M05}, Proposition 3.2$)$ $exp(\mathcal{K}_{\eta,r}(\mathbb{C}))=\mathcal{K}^{*}_{\eta,r}(\mathbb{C})$.
\end{lem}
\medskip

\begin{lem}\label{L:4}  If \  $N\in P\mathcal{K}_{\eta,r}(\mathbb{C})P^{-1}$ such that $e^{N}\in\Phi(GA(n, \mathbb{C}))$,
so $N-2ik\pi I_{n+1}\in\mathcal{F}_{n+1}$, for some $k\in\mathbb{Z}$.
\end{lem}
\medskip

\begin{proof}  Let $N'=P^{-1}NP\in\mathcal{K}_{\eta,r}(\mathbb{C})$, $M=e^{N}$ and $M'=P^{-1}MP$. We have $e^{N'}=M'$ and by Lemma ~\ref{L:4-}, $M'\in\mathcal{K}^{*}_{\eta,r}(\mathbb{C})$.
 Write $M'=\mathrm{diag}(M'_{1},\dots,M'_{r})$
 and $N'=\mathrm{diag}(N'_{1},\dots,N'_{r})$, $M'_{k}, N'_{k}\in\mathbb{T}_{n_{k}}(\mathbb{C})$, $k=1,\dots,r$. Then $e^{N'}=\mathrm{diag}(e^{N'_{1}},\dots,e^{N'_{r}})$, so
  $e^{N'_{1}}=M'_{1}$. As $1$ is the only eigenvalue of $M'_{1}$,
$N'_{1}$ has an eigenvalue $\mu\in \mathbb{C}$ such that $e^{\mu}=1$. Thus $\mu=2ik\pi$ for some $k\in\mathbb{Z}$. Therefore,
$N''=N'-2ik\pi I_{n+1}\in\mathcal{F}_{n+1}$ and satisfying $e^{N''}=e^{-2ik\pi}e^{N'}=M'$. It follows  that
$N-2ik\pi I_{n+1}=PN''P^{-1}\in P\mathcal{F}_{n+1}P^{-1}=\mathcal{F}_{n+1}$, since $P\in\Phi(GA(n, \mathbb{C}))$.
\end{proof}
\medskip

\begin{lem}\label{L:21013}$($\cite{aAh-M05}, Lemma 4.2$)$ Under above notations, one has
$exp(\mathrm{g}_{1})=G$.
\end{lem}
\medskip

As consequence, we obtain
\begin{prop}\label{p:21012}We have:\
\\ $(i)$
$\mathrm{g}_{1}= \mathrm{g}+2i\pi \mathbb{Z}I_{n+1}$.\
\\
$(ii)$ $exp(\mathrm{g})=G$.
\end{prop}
\medskip

\begin{proof} (i) By Lemma ~\ref{L:21013}, $exp(\mathrm{g}_{1})= G\subset \Phi(GA(n, \mathbb{C}))$, then by Lemma ~\ref{L:4} we have $\mathrm{g}_{1}\subset \mathrm{g}+2i\pi \mathbb{Z}I_{n+1}$.
Conversely is obvious since $\mathrm{g}+2i\pi \mathbb{Z}I_{n+1}\subset \mathcal{K}_{\eta,r}(\mathbb{C})$ and $exp(\mathrm{g}+2i\pi \mathbb{Z}I_{n+1})=exp(\mathrm{g})\subset G$.
\
\\
(ii) By Lemma ~\ref{L:21013} and (i) we have $G=exp(g_{1})=exp(\mathrm{g}+2i\pi \mathbb{Z}I_{n+1})=exp(\mathrm{g})$.
\end{proof}
\medskip

\section{ Proof of Theorem ~\ref{T:1}}
Let $\widetilde{G}$ be the group generated by $G$ and
$H=\{\lambda.I_{n+1}:\  \  \ \lambda\in \mathbb{C}^{*} \}$, and write
$\widetilde{\mathrm{g}}=exp^{-1}(\widetilde{G}\cap \mathcal{K}_{\eta,r}(\mathbb{C}))$. Then
$\widetilde{G}$ is an abelian subgroup of $GL(n+1,\mathbb{C})$. See that
$Id_{\mathbb{C}^{n+1}}\in exp^{-1}(H)\subset \widetilde{\mathrm{g}}$, so
$\widetilde{\mathrm{g}}\backslash \Psi(MA(n, \mathbb{C}))\neq\emptyset$.
\
\\
Denote by:\
\\
- $(\mathrm{g}_{1})_{v_{0}}=\{Bv_{0}\ : \ B\in \mathrm{g}_{1}\}$.\
\\
- $\widetilde{\mathrm{g}}_{v_{0}}=\{Bv_{0}\ : \ B\in \widetilde{\mathrm{g}}\}$.\
\\
\begin{prop}\label{C0C:1}$($\cite{aAh-M05},\ Theorem1.1$)$ Let  $G$ be an abelian subgroup of $GL(n+1,
\mathbb{C})$. The following are equivalent:
\
\\
(i) $G$ has a dense orbit in $\mathbb{C}^{n+1}$
\\
(ii) The orbit $G(v_{0})$ is dense in $\mathbb{C}^{n+1}$
\\
(iii) $(\mathrm{g}_{1})_{v_{0}}$ is an additive subgroup dense in
$\mathbb{C}^{n+1}$
\end{prop}
\medskip

\begin{lem}\label{LL0L:1} The following assertions are equivalent:\
\\
$(i)$ $\overline{\mathfrak{g}_{w_{0}}}=\mathbb{C}^{n}$.\
\
\\
$(ii)$ $\overline{\mathrm{g}_{v_{0}}}=\{0\}\times\mathbb{C}^{n}$.\
\\
$(iii)$ $\overline{\widetilde{\mathrm{g}}_{v_{0}}}=\mathbb{C}^{n+1}$.\

\end{lem}
\medskip

\begin{proof}$(i)\Longleftrightarrow (ii):$ For every  $f'=(B,b)\in \mathfrak{g}$, one has \begin{align*}\Psi(f')v_{0}& =\left[\begin{array}{cc}
                                                                          0 & 0 \\
                                                                          b & B
                                                                        \end{array}
\right]\left[\begin{array}{c}
               1 \\
               w_{0}
             \end{array}
\right]\\
\ & =\left[\begin{array}{c}
               0\\
               b+Bv_{0}
             \end{array}
\right].
\end{align*}
Then $\{0\}\times\mathfrak{g}_{w_{0}}=\mathrm{g}_{v_{0}}$  and the equivalence is proved.
\
\\
\\
$(ii)\Longleftrightarrow (iii):$ Firstly, remark that $\widetilde{\mathrm{g}}=\mathrm{g}_{1}+H$ and $v_{0}\in\{1\}\times\mathbb{C}^{n}$, then
$\overline{\widetilde{\mathrm{g}}_{v_{0}}}=\overline{(\mathrm{g}_{1})_{v_{0}}}+ \mathbb{C}v_{0}$. By Proposition ~\ref{p:21012}.(i),
$(\mathrm{g}_{1})_{v_{0}}=\mathrm{g}_{v_{0}}+2i\pi\mathbb{Z}v_{0}$, so $\overline{\widetilde{\mathrm{g}}_{v_{0}}}=\overline{\mathrm{g}_{v_{0}}}+ \mathbb{C}v_{0}$.
   Secondly, suppose that  $\overline{\mathrm{g}_{v_{0}}}=\{0\}\times\mathbb{C}^{n}$. Since $v_{0}\notin \{0\}\times\mathbb{C}^{n}$ and
 $Id_{\mathbb{C}^{n+1}}\in exp^{-1}(H\cap \mathcal{K}_{\eta,r}(\mathbb{C}))\subset \widetilde{\mathrm{g}}$,
 then $\mathbb{C}_{v_{0}}\subset\widetilde{\mathrm{g}}_{v_{0}}$. Therefore  $\mathbb{C}^{n+1}=\overline{\mathrm{g}_{v_{0}}}\oplus \mathbb{C}v_{0}\subset \overline{\widetilde{\mathrm{g}}_{v_{0}}}$.\
 \\
 Conversely, suppose that $\overline{\widetilde{\mathrm{g}}_{v_{0}}}=\mathbb{C}^{n+1}$. Since
 $\overline{\mathrm{g}_{v_{0}}}\subset\{0\}\times\mathbb{C}^{n}$ and $\mathbb{C}v_{0}\cap\left(\{0\}\times\mathbb{C}^{n}\right)=\{0\}$, then
 $\overline{\mathrm{g}_{v_{0}}}\oplus \mathbb{C}v_{0}=\mathbb{C}^{n+1}$,
 thus $\overline{\mathrm{g}_{v_{0}}}=\{0\}\times\mathbb{C}^{n}$.
\end{proof}
\medskip

\begin{lem}\label{LL1L:9} Let $x\in\mathbb{C}^{n}$. Then the following assertions are equivalent:\
\\
(i) $\overline{\mathcal{G}(x)}=\mathbb{C}^{n}$.\
\\
(ii) $\overline{G(1,x)}=\{1\}\times\mathbb{C}^{n}$.\
\\
(iii) $\overline{\widetilde{G}(1,x)}=\mathbb{C}^{n+1}$.\
\end{lem}
\medskip

\begin{proof} $(i)\Longleftrightarrow (ii):$ The proof is obvious, since $\{1\}\times\mathcal{G}(x)=G(1,x)$ by construction.
\
\\
$(ii)\Longleftrightarrow (iii):$ Suppose that
$\overline{\widetilde{G}(1,x)}=\mathbb{C}^{n+1}$. If  $\overline{G(1,x)}\neq \{1\}\times\mathbb{C}^{n}$,  then there
exists an open subset $O$ of $\mathbb{C}^{n}$ such that
$(\{1\}\times O)\cap G(1,x)=\emptyset$. Let $y\in O$ and $(g_{m})_{m}$ be a sequence in $\widetilde{G}$ such that
$\underset{m\to +\infty}{lim}g_{m}(1,x)=(1,y)$. Since
$\widetilde{G}$ is abelian then $g_{m}=\lambda_{m}f_{m}$, with
$f_{m}\in G$ and $\lambda_{m}\in\mathbb{C}^{*}$, thus
$\underset{m\to +\infty}{lim}\lambda_{m}=1$. Therefore,
$\underset{m\to
+\infty}{lim}f_{m}(1,x)=\underset{m\to
+\infty}{lim}\frac{1}{\lambda_{m}}g_{m}(1,x)=(1,y)$. Hence,
$(1,y)\in \overline{G(1,x)}\cap \left(\{1\}\times O\right)$, a contradiction. Conversely, if $\overline{G(1,x)}=\{1\}\times\mathbb{C}^{n}$, then \begin{align*}
\mathbb{C}^{n+1} & = \underset{\lambda\in \mathbb{C}}{\bigcup}\lambda \left(\{1\}\times\mathbb{C}^{n}\right)\\
\ & = \underset{\lambda\in \mathbb{C}}{\bigcup}\lambda \overline{G(1,x)}\\
\ & \subset \overline{\widetilde{G}(1,x)}
\end{align*}
\end{proof}

\subsection{{\bf Proof of Theorem ~\ref{T:1}}}  The proof of Theorem ~\ref{T:1} results directly  from Lemma ~\ref{LL1L:9},
 Proposition ~\ref{C0C:1} and Lemma ~\ref{LL0L:1}. $\hfill{\Box}$
\bigskip

\section{ Finitely  generated  subgroups}
\medskip

\subsection{Proof of Theorem ~\ref{T:2}} \

Let $J_{k} =\mathrm{diag}(J_{k,1},\dots, J_{k,r})$ with \ $J_{k,i}=0\in
\mathbb{T}_{n_{i}}(\mathbb{C})$ if $i\neq k$ and $J_{k,k} =
I_{n_{k}}$.
\
\\
\begin{prop}\label{p:10}$($\cite{aAh-M05}, \ Proposition 8.1$)$ Let $G$ be an abelian subgroup of
$GL(n+1,\mathbb{C})$  generated by $A_{1},\dots,A_{p}$.
Let  $B_{1},\dots,B_{p}\in \mathrm{g}$ such that  $A_{k} =
e^{B_{k}}$, $k = 1,\dots,p$ and $P\in GL(n+1,\mathbb{C})$ satisfying $P^{-1}GP\subset \mathcal{K}^{*}_{\eta,r}(\mathbb{C})$. Then:
 $$\mathrm{g}_{1} =
\underset{k=1}{\overset{p}{\sum}}\mathbb{Z}B_{k}+2i\pi\underset{k=1}{\overset{r}{\sum}}\mathbb{Z}PJ_{k}P^{-1} \ \  \mathrm{and} \ \
 \ (\mathrm{g_{1}})_{v_{0}} = \underset{k=1}{\overset{p}{\sum}}\mathbb{Z}B_{k}v_{0} + \underset{k=1}{\overset{r}{\sum}} 2i\pi \mathbb{Z}Pe^{(k)}.$$
 \end{prop}
\bigskip

\begin{prop} \label{p:11} $($Under notations of Proposition ~\ref{p:2}$)$ Let  $\mathcal{G}$ be an abelian subgroup of  $GA(n, \mathbb{C})$ generated by
 $f_{1},\dots,f_{p}$ and let  $f'_{1},\dots,f'_{p}\in
\mathfrak{g}$ such that $\Phi(f_{k}) = e^{\Psi(f'_{k})}$, $k =
1,..,p$. Then:
 $$\mathfrak{g}_{w_{0}} = \underset{k=1}{\overset{p}{\sum}}\mathbb{Z}f'_{k}(w_{0}) + \underset{k=1}{\overset{r}{\sum}} 2i\pi \mathbb{Z}(p_{2}\circ Pe^{(k)}).$$
 \end{prop}
\medskip

\begin{proof} Let  $G=\Phi(\mathcal{G})$.  Then  $G$ is
generated by  $\Phi(f_{1}),\dots,\Phi(f_{p})$.\ By proposition ~\ref{p:10} we have
$$\mathrm{g_{1}} =
\underset{k=1}{\overset{p}{\sum}}\mathbb{Z}\Psi(f'_{k}) +
\underset{k=1}{\overset{r}{\sum}} 2i\pi \mathbb{Z}PJ^{(k)}P^{-1}.$$

Since $P\in \Phi(GA(n, \mathbb{C}))$ then $PJ^{(1)}P^{-1}\notin \mathcal{F}_{n+1}$
and $PJ^{(k)}P^{-1}\in \mathcal{F}_{n+1}$ for every $k=2,\dots,r$.
As  $\mathrm{g}=\mathrm{g}_{1}\cap \mathcal{F}_{n+1}$, then
$$\mathrm{g} =\left\{\begin{array}{c}
                       \underset{k=1}{\overset{p}{\sum}}\mathbb{Z}\Psi(f'_{k}) +
\underset{k=2}{\overset{r}{\sum}} 2i\pi \mathbb{Z}PJ^{(k)}P^{-1}, \ if\ r\geq2 \\
                       \underset{k=1}{\overset{p}{\sum}}\mathbb{Z}\Psi(f'_{k}), \ if\ r=1 \ \ \ \ \ \ \ \ \ \ \ \ \ \ \ \ \ \ \ \ \
                     \end{array}
\right.$$
By construction, one has
$\mathfrak{g}_{w_{0}}=p_{2}\left(\mathrm{g}_{v_{0}}\right)$, $v_{0}=Pu_{0}$, $J^{(k)}u_{0}=e^{(k)}$ and $p_{2}(\Psi(f'_{k})v_{0})=f'_{k}(w_{0})$.
Then $$\mathfrak{g}_{w_{0}}=\left\{\begin{array}{c}
                                     \underset{k=1}{\overset{p}{\sum}}\mathbb{Z}f'_{k}(w_{0}) + \underset{k=2}{\overset{r}{\sum}} 2i\pi \mathbb{Z}(p_{2}\circ Pe^{(k)}),\ if\ r\geq2 \\
                                     \underset{k=1}{\overset{p}{\sum}}\mathbb{Z}f'_{k}(w_{0}) , \ if\ r=1\ \ \ \ \ \ \ \ \  \ \ \ \ \ \ \ \ \ \  \ \ \ \ \ \ \ \ \ \ \
                                   \end{array}
\right.\ \ \ \ (3)$$
The proof is completed.
\end{proof}

Recall the following Proposition which was proven in \cite{mW}:

\begin{prop}\label{p:12}$(cf.$ \cite{mW}, $page$ \ $35)$. Let $F = \mathbb{Z}u_{1}+\dots+\mathbb{Z}u_{p}$ with $u_{k } =
(u_{k,1}, \dots,u_{k,n})\in\mathbb{C}^{n}$ and $u_{k,i} =
Re(u_{k,i})+i Im(u_{k,i})$, $k = 1,\dots, p$, $i = 1,\dots,n$. Then
$F$ is dense in $\mathbb{C}^{n}$ if and only if for every
$(s_{1},\dots,s_{p})\in
 \mathbb{Z}^{p}\backslash\{0\}$ :
$$rank\left[\begin{array}{cccc }
 Re(u_{1,1 }) &\dots &\dots & Re(u_{p,1 }) \\
 \vdots &\vdots &\vdots &\vdots \\
 Re(u_{1, n }) &\dots &\dots& Re(u_{p,n }) \\
 Im(u_{1,1 }) &\dots &\dots & Im(u_{p, 1 }) \\
 \vdots &\vdots&\vdots &\vdots \\
 Im(u_{1, n }) &\dots &\dots & Im(u_{p, n }) \\
 s_{1 } &\dots&\dots & s_{p }
 \end{array}\right] =\ 2n+1.$$
\end{prop}
\
\\
\\
\emph{Proof of Theorem ~\ref{T:2}:} This follows directly from Theorem
~\ref{T:1}, Propositions ~\ref{p:11} and ~\ref{p:12}.
\bigskip
\
\\
\\
\subsection{Proof of Corollaries ~\ref{C:9} and ~\ref{C:10}} \
\\
\emph{Proof of Corollary ~\ref{C:9}:} We show first that if $F = \mathbb{Z}u_{1}+\dots +\mathbb{Z}u_{m}$,
$u_{k}\in \mathbb{C}^{n}$ with $m \leq 2n$, then $F$ can not be dense: Write
$u_{k} \in \mathbb{C}^{n}$, $u_{k} = Re(u_{k}) + iIm(u_{k})$ and
$v_{k} = [Re(u_{k}); Im(u_{k}); s_{k}]^{T}\in \mathbb{R}^{2n+1}$,
 $1\leq k \leq m$. Since $m \leq 2n$, it follows that $rank(v_{1},\dots, v_{m}) \leq 2n$,
and so $F$ is not dense in $\mathbb{C}^{n}$ by Proposition ~\ref{p:12}.\\
Now, by applying Theorem ~\ref{T:2} and the fact that $m = p + r-1 \leq 2n$ (since $r\leq n+1$) and by (1), the
Corollary ~\ref{C:9} follows.  $\hfill{\Box}$
\
\\
\\
\emph{Proof of Corollary ~\ref{C:10}:} Since $p\leq n$ and $r\leq n+1$ then $p + r-1\leq 2n$.
Corollary ~\ref{C:10}  follows from Corollary ~\ref{C:9}. $\hfill{\Box}$

\medskip
\section{ Example}
\bigskip

\begin{exe} \label{exe:2} Let  $\mathcal{G}$\ the group generated by
$f_{1}=(A_{1},a_{1})$, $f_{2}=(A_{2},a_{2})$,
$f_{3}=(A_{3},a_{3})$ and $f_{4}=(A_{4},a_{4})$, where: $$A_{1}=I_{2}, \ \ a_{1}=(1+i, 0),\ \ \ \ A_{2}=\mathrm{diag}(1, e^{-2+i}),
 \ \ a_{2}=\left(0, \ 0\right).$$
$$
A_{3} =\mathrm{diag}\left(1,\
e^{\frac{-\sqrt{2}}{\pi}+i\left(\frac{\sqrt{2}}{2\pi}-\frac{\sqrt{7}}{2}\right)}\right), \ \ \ \ \
 a_{3} =\left(\frac{-\sqrt{3}}{2\pi}+i\left(\frac{\sqrt{5}}{2}-\frac{\sqrt{3}}{2\pi}\right),0\right),$$
$$A_{4}=I_{2}\ \ \ \  \mathrm{and} \ \ \ \ \ \ \ a_{4}=(2i\pi,\ 0).$$

Then every orbit in  $V=\mathbb{C}\times \mathbb{C}^{*}$ is dense in  $\mathbb{C}^{2}$.
\end{exe}

\bigskip

\begin{proof} Denote by  $G=\Phi(\mathcal{G})$.  Then  $G$ generated by $$\Phi(f_{1})=\left[\begin{array}{ccc}
                                                                                                    1 & 0 & 0 \\
                                                                                                    1+i & 1 & 0\\
                                                                                                    0 & 0 & 1
                                                                                                  \end{array}
\right], \Phi(f_{2})=\left[\begin{array}{ccc}
 1 & 0 & 0 \\
 0 & 1 & 0\\
 0 & 0 & e^{-2+i}
 \end{array}
\right],$$
\ $$\Phi(f_{3})=\left[\begin{array}{ccc}
                               1 & 0 & 0 \\
                               \frac{-\sqrt{3}}{2\pi}+i\left(\frac{\sqrt{5}}{2}-\frac{\sqrt{3}}{2\pi}\right) & 1 & 0 \\
                               0 & 0 & e^{\frac{-\sqrt{2}}{\pi}+i\left(\frac{\sqrt{2}}{2\pi}-\frac{\sqrt{7}}{2}\right)}
                              \end{array}
\right],$$ and $$\Phi(f_{4})=\left[\begin{array}{ccc}
                               1 & 0 & 0 \\
                              2i\pi & 1 & 0 \\
                               0 & 0 & 1
                             \end{array}
\right].$$

Let $f'_{1}=(B_{1},b_{1})$, $f'_{2}=(B_{2},b_{2})$ and $f'_{2}=(B_{2},b_{2})$  such that $e^{\Psi(f'_{k})}=A'_{k}$, $k=1,2,3,4$. We have
 $$B_{1}=\mathrm{diag}(0,\ 0),  b_{1}=(1+i,\ 0)\ \ ,$$
$$B_{2}=\mathrm{diag}(0, \ -2+i),\ \ \ \ b_{2}=(0,\ 0),$$
$$B_{3}=\mathrm{diag}\left(0,\ \ \
   \frac{-\sqrt{2}}{\pi}+i\left(\frac{\sqrt{2}}{2\pi}-\frac{\sqrt{7}}{2}\right)\right),\ \  \ \ b_{3}=\left(\frac{-\sqrt{3}}{2\pi}+i\left(\frac{\sqrt{5}}{2}-\frac{\sqrt{3}}{2\pi}\right),\ 0\right),$$
    $$B_{4}=\mathrm{diag}(0,0) \ \ \ \ \  \mathrm{and} \ \ \ \ \ \ b_{4}=(2i\pi, \ 0).$$
\
\\
Here, we have:\
\\
- $G$ is an abelian subgroup of $\mathcal{K}_{(2,1),2}^{*}(\mathbb{C})$.\
\\
- $P=I_{3}$, $r=2$, $U=\mathbb{C}^{*}\times\mathbb{C}\times\mathbb{C}^{*}$, $u_{0}=(1,0,1)$, $e^{(1)}=(1,0,0)$ and $e^{(2)}=(0,0,1)$\
\\
- $V=\mathbb{C}\times \mathbb{C}^{*}$, $w_{0}=(0,1)$.
\
\\
\\
In the other hand, for every $(s_{1}, s_{2}, s_{3}, s_{4},
t_{2})\in\mathbb{Z}^{5}\backslash\{0\}$, one has the determinant:

\begin{align*}
\Delta & = \left|\begin{array}{ccccc }
 Re(B_{1}w_{0}+b_{1}) & Re(B_{2}w_{0}+b_{2}) & Re(B_{3}w_{0}+b_{3}) & Re(B_{4}w_{0}+b_{4}) & 0 \\
 \\
 Im(B_{1}w_{0}+b_{1}) & Im(B_{2}w_{0}+b_{2}) & Im(B_{3}w_{0}+b_{3}) & Im(B_{4}w_{0}+b_{4}) & 2\pi e^{(2)} \\
 \\
  s_{1 } & s_{2 } & s_{3 } & s_{4 } & t_{2 }
  \end{array}\right|  \\
\ & =\left|\begin{array}{ccccc }
 -\frac{\sqrt{3}}{2\pi } & 1 & 0 & 0 & 0 \\
  -\frac{\sqrt{2}}{\pi } & 0 & -2 & 0 & 0 \\
  \frac{\sqrt{5}}{2}-\frac{\sqrt{3}}{2\pi } & 1 & 0 & 2\pi & 0 \\
  \frac{\sqrt{2}}{2\pi}-\frac{\sqrt{7}}{2 } & 0 & 1 & 0 & 2\pi \\
  s_{1 } & s_{2 } & s_{3 } & s_{4 } & t_{2 }
  \end{array}\right|\\
  \ & =-2\pi\left((4s_{1})\pi-(2s_{3})\sqrt{2}+(2s_{2})\sqrt{3}+s_{4}\sqrt{5}+t_{2}\sqrt{7}\right).
  \end{align*}

Since $\pi$, $\sqrt{2}$, $\sqrt{3}$, $\sqrt{5}$ and $\sqrt{7}$ \ \
are rationally independent then $\Delta\neq 0$ for every
 $(s_{1}, s_{2}, s_{3}, t_{1}, t_{2})\in\mathbb{Z}^{5}\backslash\{0\}$. It follows that: $$rank\left[\begin{array}{ccccc }
 -\frac{\sqrt{3}}{2\pi } & 1 & 0 & 0 & 0 \\
  -\frac{\sqrt{2}}{\pi } & 0 & -2 & 0 & 0 \\
  \frac{\sqrt{5}}{2}-\frac{\sqrt{3}}{2\pi } & 1 & 0 & 2\pi & 0 \\
  \frac{\sqrt{2}}{2\pi}-\frac{\sqrt{7}}{2 } & 0 & 1 & 0 & 2\pi \\
  s_{1 } & s_{2 } & s_{3 } & s_{4 } & t_{2}
  \end{array}\right]\ =\ 5$$
  and by Theorem ~\ref{T:2}, $\mathcal{G}$ has a dense orbit and every orbit of $V$ is dense in $\mathbb{C}^{2}$.
\end{proof}

\bigskip

\bibliographystyle{amsplain}

\begin{thebibliography}{10}

\bibitem{aAh-M05} A. Ayadi and H. Marzougui, \emph{Dense orbits for abelian subgroups of
GL(n, C)}, Foliations 2005: 47--69. World
Scientific,Hackensack,NJ,2006.

\bibitem{bm} F. Bayart, E. Matheron, Dynamics of Linear Operators, Cambridge Tracts in Math., 179, Cambridge University Press, 2009.

\bibitem{Ja} M. Javaheri, \emph{A generalization of Dirichlet approximation theorem
for the affine actions on real line}, Journal of Number Theory \textbf{128} (2008) 1146--1156.


\bibitem{Ku} R.S. Kulkarni, \emph{Dynamics of linear and affine maps}, arXiv:0709.1353v1, 2007.

\bibitem{Bu} T. V. Budnytska, \emph{Classification of topologically conjugate affine mapping}, Ukrainian Mathematical Journal, \textbf{61}, No. 1, (2009) 164--170.

\bibitem{be} V. Bergelson, M. Misiurewicz and S. Senti,  \emph{Affine actions of a free semigroup on the real line}
Ergod. Th. and Dynam. Sys. \textbf{26}, (2006), 1285--1305.

\bibitem{mW} Waldschmidt.M, \emph{Topologie des points rationnels},
Cours de troisi\`eme Cycle, Universit\'e P. et M. Curie (Paris
VI), 1994/95.
\end{thebibliography}
\vskip 0,4 cm
%%% ----------------------------------------------------------------------

\end{document}